\documentclass[11pt,reqno]{amsart}
\usepackage{verbatim}
\usepackage{amssymb}
\usepackage{amsmath}
\usepackage{graphicx}
\usepackage{appendix}
\usepackage{color}
\usepackage{amsthm}
\usepackage{tikz}
\usepackage{mathrsfs}
\usetikzlibrary{arrows,positioning,decorations.pathmorphing,decorations.markings}

\def\t{\tau}

\def\w{\omega}

\def\=>{\Longrightarrow}

\def\to{\longrightarrow}

\def\^{\wedge}
\def\+{\dagger}

\def\dd[#1,#2]{\frac{d#1}{d#2}}
\def\del[#1,#2]{\frac{\partial #1}{\partial #2}}
\def\over[#1]{\overline{#1}}
\def\vec[#1]{\overrightarrow{#1}}

\tikzset{->-/.style={decoration={
  markings,
  mark=at position .5 with {\arrow{latex}}},postaction={decorate}}}
\tikzset{
    >=latex
    }
\begin{document}

%% --------------------------------------------------------------------------
%% Please use the environment "talk" for each abstract.
%% It has three obligatory and one optional argument. The syntax is:
%% -----------------------
%% \begin{talk}[coauthors]{Name of the speaker}{Title of the talk}{Author Sorting Index}
%%      .....
%% \end{talk}
%% -----------------------
%% The names of coauthors will appear in form of "(joint work with ...)"
%%
%% The Author Sorting Index should be given as the last and first name of the speaker,
%% separated by a comma. If for example the name of the speaker is "John Smith", then
%% the correct Author Sorting Index is "Smith, John".
%% Any special characters (like accents, German umlaute, etc.) should be replaced by
%% their "non-special" version, eg replace \"a by a, \'a by a, etc.
%%
%% Please use the standard thebibliography environment to include
%% your references, and try to use labels for the bibitems, which
%% are uniquely assigned to you in order to avoid conflicts with other authors.
%% You can achieve unique labels by using our on initials before every label.
%% -------------------------------------------------------------------------------

\title{Super-Teichm\"uller spaces and related structures}

\author{Anton M. Zeitlin}
\begin{abstract}
This short note provides an overview of some theorems and conjectures obtained by the author and his collaborators. It is an extended abstract for the Oberwolfach workshop ``New Trends in Teichm\"uller Theory and Mapping Class Groups", 2 September - 8 September 2018.
\end{abstract}

\date{\today}

\numberwithin{equation}{section}

\maketitle

\section { Introduction.} The study of the superstring theory has drawn the attention to the very important generalizations of Riemann surfaces, known as super-Riemann surfaces (SRS) that could be viewed as $(1|N)$-dimensional complex supermanifolds with extra structures \cite{CR}, \cite{schwarz}. Moduli spaces of $N=1$ SRS are of special importance (see e.g. \cite{witten} for a review).  

One of the ways to look at the corresponding Teichm\"uller spaces $S\mathcal{T}(F)$ (here $F$ is the underlying Riemann surface of genus $g$ with $s$ punctures),  is the view through the prism of the higher Teichm\"uller theory. They are defined in a manner analogous to that for the standard pure 
even case, ${S\mathcal{T}(F)=\mathrm{Hom}'(\pi_1(F)\to G)/G}$,  where instead of $PSL(2,\mathbb{R})$, $G$ is a {\it supergroup} $OSp(1|2)$.  Here $\pi_1(F)$ is the fundamental group of the underlying Riemann surface with punctures, and $\mathrm{Hom}'$ stands for the homomorphisms that map the elements of $\pi_1(F)$, corresponding to small loops around the punctures, to parabolic elements of $OSp(1|2)$, which means that their natural projections to $PSL(2,\mathbb{R})$ are parabolic elements. 

 The image of the fundamental group under ${\rm Hom}'$ produces a generalization of the standard Fuchsian group $\Gamma$, which acts on a super-analogue of the upper half-plane $H^+$ producing $N=1$ super-Riemann surfaces as a factor $H^+/\Gamma$ \cite{CR}. It is necessary to build the super analogues of known objects  in Teichm\"uller theory for the successful study of such spaces. In this note I will give an overview of the results from \cite{PZ}, \cite{IPZ}, \cite{IPZ2} concerning the generalization of the Penner coordinates \cite{penner}, \cite{pbook} for the super-Teichm\"uller space. 
 
The Penner coordinates are the coordinates on $\mathbb{R}^s_+$-bundle 
$\tilde{\mathcal{T}}(F)=\mathbb{R}^s_+\times\mathcal{T}(F)$ over the  Teichm\"uller space of $\mathcal{T}(F)$ of $s$-punctured surfaces with negative Euler characteristics. 

The construction is based on  the {\it ideal triangulation} of $F$  (i.e.  vertices of triangualtion are the punctures) and the assigninment of a positive number to every edge of the triangulation. An important feature of these coordinates is that under the elementary changes of triangulation, known as Whitehead moves, or {\it flips} generating the mapping class group,  
the change of coordinates is rational, described by the so-called {\it Ptolemy relations}. Therefore,  it is  making the mapping class group action rational. 

The difficulty in constructing an analogue of these coordinates for $S\tilde{\mathcal{T}}(F)=\mathbb{R}^s_+\times S{\mathcal{T}}(F)$ is that $S{\mathcal{T}}(F)$ has many connected components enumerated by spin structures. Thus, to proceed further it is necessary to have
a suitable combinatorial description of the spin structures.

\section {Fatgraphs and spin structures}  Consider the trivalent {\it fatgraph} $\tau$, corresponding to an $s$-punctured $(s>0)$ Riemann surface $F$ of the negative Euler characteristic (i.e. a graph with trivalent vertices), which is homotopically equivalent to $F$, with cyclic orderings on half-edges for every vertex \cite{pbook} induced by the orientation of the surface. 

There is a one-to-one correspondence between the ideal triangulations and trivalent fatgraphs. Let $\w$ be an orientation on the edges of $\t$. As in \cite{PZ}, we define a \emph{fatgraph reflection} at a vertex $v$ of $(\t,\w)$ as a reversal of the orientations of $\w$ on every edge of $\t$ incident to $v$. Let us define the $\mathcal{O}(\t)$ to be the equivalence classes of orientations on a trivalent fatgraph $\tau$ of $F$, where $\omega_1\sim \omega_2$ if and only if $\omega_1$ and $\omega_2$ differ by a finite number of fatgraph reflections. In \cite{PZ},\cite{IPZ}, we identified such classes of orientations on fatgraphs with the spin structures on $F$. The paths corresponding to the boundary cycles on the fatgraph (i.e. the punctures of $F$) are divided into two classes depending on the parity of number $k$ -- the number of edges with orientation opposite to the canonical orientation of $\gamma$. 

\begin{figure}[h!]
\begin{center}
\begin{tikzpicture}[ scale=0.6, ultra thick, baseline=1cm]
\draw (0,0)--(210:1) node[above] at (210:0.7){$\epsilon_2$};
\draw (0,0)--(330:1) node[above] at (330:0.7){$\epsilon_4$};
    % draw the connecting line
    \draw[ 
 	ultra thick,
        decoration={markings, mark=at position 0.6 with {\arrow{>}}},
        postaction={decorate}
        ]
        (0,0) -- (0,2);
\draw[yshift=2cm] (0,0)--(30:1) node[below] at (30:0.7){$\epsilon_3$};
\draw[yshift=2cm] (0,0)--(150:1) node[below] at (150:0.7){$\epsilon_1$};
\end{tikzpicture}
\begin{tikzpicture}[scale=0.6, baseline]
\draw[->, thick](0,0)--(2,0);
\node[above] at (2,0) {};
\node at (-1,0){};
\node at (3,0){};
\end{tikzpicture}
\begin{tikzpicture}[scale=0.7, ultra thick, baseline]
\draw (0,0)--(120:1) node[above] at (100:0.3){$\epsilon_1$};
\draw (0,0)--(240:1) node[below] at (260:0.3){$\epsilon_2$};
    % draw the connecting line
    \draw[ 
        decoration={markings, mark=at position 0.5 with {\arrow{<}}},
        postaction={decorate}
        ]
        (0,0) -- (2,0);
\draw[xshift=2cm] (0,0)--(60:1) node[above] at (100:0.4){$-\epsilon_3$};
\draw [xshift=2cm](0,0)--(-60:1) node[below] at (-80:0.3){$\epsilon_4$};
\end{tikzpicture}
\end{center}
\caption{Spin graph evolution in the generic situation}
\label{flipgraphint}
\end{figure}
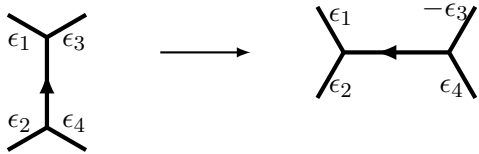

The punctures are called Ramond (R) when $k$ is even, and Neveu-Schwarz (NS) \cite{witten} when $k$ is odd. In \cite{PZ}, we have also proved that under the flip transformations the orientations change  in the generic situation as in Figure \ref{flipgraphint}, 
where $\epsilon_i$ stand for orientations on edges, and extra minus sign stands for the orientation reversal.
% In fact there is another way of thinking about the spin structures, using graph connections, which we will need later. 
%If $G$ be a group, a \emph{$G$-graph connection} on $\t$ is the assignment $g_e\in G$ to each oriented edge $e$ of $\t$ so that $g_{\over[e]}=g_e\inv$ if $\over[e]$ is the opposite orientation to $e$. Two assignments $\{g_e\},\{g_e'\}$ are equivalent iff there are $t_v\in G$ for each vertex $v$ of $\t$ such that $g_e'=t_v g_e t_w\inv$ for each oriented edge $e\in\t_1$ with initial point $v$ and terminal point $w$.

%Therefore, the space of spin structures on $F$ is identified with $\mathbb{Z}_2$-graph connections on a given fathraph $\tau$ of $F$.

\section{Main Result.}
 In \cite{PZ}, we have proved the following Theorem. 
 
\noindent {\bf Theorem. }{\it 
i) The components of $S\tilde{T}(F)$ are determined by the space of spin structures on $F$. For each component $C$ of $S\tilde{T}(F)$, there are global affine coordinates on $C$ given by assigning to a triangulation $\Delta$ of $F$, 
\begin{itemize}
\item one even coordinate called $\l$-length for each edge; 
\item one odd coordinate called $\mu$-invariant for each triangle, taken modulo an overall change of sign.
\end{itemize} In particular we have a real-analytic homeomorphism: 
$${C\to \mathbb{R}_{>0}^{6g-6+3s|4g-4+2s}/\mathbb{Z}_2.}$$
ii) The super Ptolemy transformations \cite{PZ} provide the analytic relations between coordinates assigned to different choice of triangulation $\Delta'$ of $F$, namely upon flip transformation. Explicitly (see Figure \ref{ptolemy}),  when all $a,b,c,d$ are different edges of the triangulations of $F$, the Ptolemy transformations are as follows:

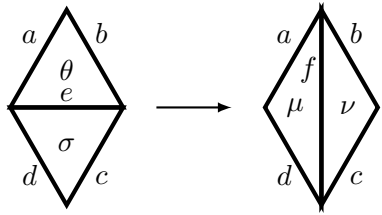
\begin{figure}[h!]

\centering

\begin{tikzpicture}[scale=0.5, baseline,ultra thick]

\draw (0,0)--(3,0)--(60:3)--cycle;

\draw (0,0)--(3,0)--(-60:3)--cycle;

\draw node[above] at (70:1.5){$a$};

\draw node[above] at (30:2.8){$b$};

\draw node[below] at (-30:2.8){$c$};

\draw node[below=-0.1] at (-70:1.5){$d$};

\draw node[above] at (1.5,-0.15){$e$};

\draw node[left] at (0,0) {};

\draw node[above] at (60:3) {};

\draw node[right] at (3,0) {};

\draw node[below] at (-60:3) {};

\draw node at (1.5,1){$\theta$};

\draw node at (1.5,-1){$\sigma$};

\end{tikzpicture}
\begin{tikzpicture}[baseline]

\draw[->, thick](0,0)--(1,0);

\node[above]  at (0.5,0) {};

\end{tikzpicture}
\begin{tikzpicture}[scale=0.5, baseline,ultra thick]

\draw (0,0)--(60:3)--(-60:3)--cycle;

\draw (3,0)--(60:3)--(-60:3)--cycle;

\draw node[above] at (70:1.5){$a$};

\draw node[above] at (30:2.8){$b$};

\draw node[below] at (-30:2.8){$c$};

\draw node[below=-0.1] at (-70:1.5){$d$};

\draw node[left] at (1.7,1){$f$};

\draw node[left] at (0,0) {};

\draw node[above] at (60:3) {};

\draw node[right] at (3,0) {};

\draw node[below] at (-60:3) {};

\draw node at (0.8,0){$\mu$};

\draw node at (2.2,0){$\nu$};

\end{tikzpicture}
\caption{Generic flip transformation}
\label{ptolemy}
\end{figure}
$$
ef=(ac+bd)\Big(1+\frac{\sigma\theta\sqrt{\chi}}{1+\chi}\Big),\quad \nu=\frac{\sigma-\theta\sqrt{\chi}}{\sqrt{1+\chi}},\quad
\mu=\frac{\theta+\sigma\sqrt{\chi}}{\sqrt{1+\chi}},
$$
where $\chi=\frac{ac}{bd}$,  so that the evolution of arrows is as in  Figure \ref{flipgraphint}.}

To every edge $e$ (see the picture above) we can associate the {\it shear coordinate} $z_e=\log (\frac{ac}{bd})$. These parameters satisfy a linear relation for every puncture, and together with odd variables they form a set of coordinates on the $ST(F)$, thus producing the removal of the decoration, as it was in the pure even case. Here we mention that there is a physically and algebro-geometrically interesting refinement of $ST(F)$ studied in \cite{IPZ2}, corresponding to the removal of certain odd degrees of freedom associated with R-punctures.

\section{ $N=2$ Super-Teichm\"uller space and beyond.} Replacing $OSp(1|2)$ in the definition of $ST(F)$ by $OSP(2|2)$, one obtains the super-Teichm\"uller space of punctured $N=2$ SRS.  It has been investigated in \cite{IPZ}, and the analogue of Penner coordinates was constructed there. Unlike $N=1$ case, the resulting $N=2$ SRS, obtained by uniformization, correspond to a certain subspace in the moduli space of $N=2$ SRS \cite{schwarz}. 
We also mention that according to 
\cite{schwarz} $N=2$ SRS are in one-to-one correspondence with $(1|1)$-dimensional supermanifolds. 
  
An important problem \cite{sz} will be to see explicitly how to glue $N=1$ and $N=2$ super-Riemann surfaces using the fatgraph data following the analogue with the Strebel theory. 

Another important task is to understand the (complexified) version of the results of  \cite{PZ}, \cite{IPZ} in the  context of spectral networks and the abelianization construction of Gaiotto, Moore and Neitzke \cite{Ga}. 
In the super case it looks that only quasi-abelianization seems to work: one should be able to describe our constructions via  the moduli space of $GL(1|1)$ local systems. 

Finally, the super-Ptolemy transformations from the Theorem above, discovered in \cite{PZ} (see also  \cite{IPZ} for $N=2$ case) should lead to new interesting generalizations of cluster algebras.

\end{document}